\documentclass[11pt]{article}
\usepackage{amsmath, amssymb, amsthm}
\usepackage{verbatim}
\usepackage{multicol}
\usepackage{enumerate}
\usepackage{comment}

\usepackage{pgf}
\usepackage{tikz}
\usepackage[FIGTOPCAP]{subfigure}
\usetikzlibrary{math,calc,intersections}
\usetikzlibrary{positioning,arrows,shapes,decorations.markings,decorations.pathreplacing,matrix,patterns}
\tikzstyle{vertex}=[circle,draw=black,fill=black,inner sep=0,minimum size=3pt,text=white,font=\footnotesize]
\usepackage{cleveref}

\date{}
\title{A note on visible islands}

\author{Sophie Leuchtner\thanks{ University of California at San Diego, La Jolla, CA, 92093 USA.  Supported by a UCSD Undergraduate Summer Research Award.
 Email:
{\tt sleuchtn@ucsd.edu.}} \and Carlos M. Nicol\'as\thanks{Department of Mathematics, Virginia Tech., Blacksburg, VA, 24061 USA.  Email: {\tt cnicolas@vt.edu}} \and Andrew Suk\thanks{Department of Mathematics, University of California at San Diego, La Jolla, CA, 92093 USA. Supported an NSF CAREER award, NSF award DMS-1952786, and an Alfred Sloan Fellowship. Email: {\tt asuk@ucsd.edu}} }

\theoremstyle{plain}
\newtheorem{theorem}{Theorem}

\newtheorem{lemma}[theorem]{Lemma}
\newtheorem{conjecture}[theorem]{Conjecture}
\newtheorem{problem}[theorem]{Problem}

\theoremstyle{definition}

\newcommand{\conv}{\textnormal{conv}}

\begin{document}
\maketitle

\sloppy

\begin{abstract}
 Given a finite point set $P$ in the plane, a subset $S \subseteq P$ is called an \emph{island} in $P$ if $\conv(S) \cap P = S$.  We say that $S\subset P$ is a \emph{visible island} if the points in $S$ are pairwise visible and $S$ is an island in $P$.   The famous Big-line Big-clique Conjecture states that for any $k \geq 3$ and $\ell \geq 4$, there is an integer $n = n(k,\ell)$, such that every finite set of at least $n$ points in the plane contains $\ell$ collinear points or \emph{k} pairwise visible points.  In this paper, we show that this conjecture is false for visible islands, by replacing each point in a Horton set by a triple of collinear points.  Hence, there are arbitrarily large finite point sets in the plane with no 4 collinear members and no visible island of size $13$.

\end{abstract}

\section{Introduction}

Given a finite point set $P$ in the plane, two points $p,q \in P$ are \emph{visible} in $P$ if no other point in $P$ lies in the interior of the segment $\overline{pq}$.  The famous Big-line Big-clique Conjecture, introduced by K\'ara, P\'or, and Wood \cite{KPW}, states that for any $k$ $\geq$ 3 and $\ell$ $\geq$ 3, there is an integer $n = n(k,\ell)$, such that every finite set of at least $n$ points in the plane contains either $\ell$ collinear points or $k$ pairwise visible points. Clearly, the conjecture holds when $k = 3$ or $\ell = 3$.  K\'ara, P\'or, and Wood showed that the conjecture holds when $k = 4$ and $\ell \geq 4$, and Abel et al.~\cite{abel} verified the conjecture for $k = 5$ and $\ell \geq 4$.  The conjecture remains open for all $k\geq 6$ and $\ell \geq 4$.  See \cite{PW,mat2} for more related results.

A natural approach to the Big-line Big-clique Conjecture is to find \emph{holes} in planar point sets.  Given a finite point set $P$ in the plane, a $k$-subset $Q\subset P$ is called a \emph{$k$-hole} in $P$ if $Q$ is in convex position and $\conv(Q)\cap P = Q$.  Clearly, if $Q$ is a $k$-hole in $P$, then $Q$ consists of $k$ pairwise visible points in $P$.  Does every sufficiently large finite point set $P$ in the plane contain $\ell$ collinear points or a $k$-hole?  In \cite{abel}, Abel et al.~proved this to be true when $k = 5$, and conjectured it to be true when $k = 6$.  However, a famous construction due to Horton \cite{H} shows that this is false for $k \geq 7$ (See also Chapter 3 in \cite{mat}).  In this paper, we study a relaxed version of this question by replacing holes with \emph{visible islands}.

   Given a finite point set $P$ in the plane, a subset $S \subseteq P$ is called an \emph{island} in $P$ if $\conv(S) \cap P = S$.  We say that $S\subset P$ is a \emph{visible island} if the points in $S$ are pairwise visible and $S$ is an island in $P$.

\begin{problem}\label{problem}
Given integers $k,\ell \geq 4$, is there an integer $n = n(k,\ell)$ such that every $n$-element planar point set $P$ contains either $\ell$ collinear points or a visible island of size $k$?

\end{problem}

For $\ell \leq 3$, clearly we have $n(k,\ell) = k$.  For $\ell \geq 4$ and $k \leq 5$, $n(k,\ell)$ exists by the result of Abel et al. \cite{abel} stated above.  Our main result shows that by modifying Horton's construction \cite{H}, $n(k,\ell)$ does not exist for $\ell = 4$ and $k \geq 13.$

    \begin{theorem}\label{main}
     There exist arbitrarily large, finite point sets in the plane with no 4 collinear points and no visible island of size 13.  
    \end{theorem}

    \noindent When $\ell \geq 4$ and $6 \leq k \leq 12$, Problem \ref{problem} remains open. 
    
    \section{Proof of Theorem \ref{main}}
    
 Let us begin by recalling the definition of Horton sets. Given finite point sets $P$ and $Q$ in the plane, We say that $P$ is \emph{high above} $Q$ (or, equivalently, $Q$ is \emph{deep below} $P$) if each line determined by two points of $P$ lies above all the points of $Q$, and each line determined by two points of $Q$ lies below all of the points of $P$.   \noindent Finally, given a point $p$ in the plane, we denote $x(p)$ to be the $x$-coordinate of $p$.  Throughout the proof, we will only consider point sets whose members have distinct $x$-coordinates.

     For $n \geq 0$, a \emph{Horton set} $H_n$ is a set of $2^n$ points in the plane with no three collinear members, defined recursively as follows.   Set $H_0$ to be a single point in the plane.  Having constructed $H_{n-1} = \{p_1,p_2,\ldots, p_{2^{n-1}}\}$, whose elements are ordered by increasing $x$-coordinate, set
    
    $$H^{(1)}_{n -1 } =  \{p_1,p_2,\ldots, p_{2^{n-1}}\},$$ 
    
    $$H^{(2)}_{n-1 } =  H^{(1)}_{n-1} + (\varepsilon, K),$$

 \noindent    where $K$ is a sufficiently large number such that $H^{(1)}_{n - 1}$ lies deep below  $H^{(2)}_{n - 1}$.  Likewise, we set $\varepsilon > 0$ to be sufficiently small such that for each $i$, 
    
    $$x(p_i) < x(p_i) + \varepsilon <   x(p_{i + 1}).$$
    
  \noindent  Then we set $H_{n } = H^{(1)}_{n -1}\cup H^{(2)}_{n-1 }.$  It is known that $H_n$ has the following property (see Chapter 3 in \cite{mat}).
 
\begin{lemma}[\cite{H,mat}]\label{lem} If $S\subset H_n$ such that $|S| = 7$, then the interior of $\conv(S)$ contains a point from $H_n$.

\end{lemma}

For each $p_i\in H_n$, replace $p_i$ with three collinear points $q_i,u_i,v_i$
that are very close together, while avoiding the creation of four
collinear points.  The points $q_i,u_i,v_i$ are called \emph{triplets} of each other and $p_i$ is the \emph{parent} of them.  Let $\hat{H}_n$ be the resulting set. Then $\hat{H}_n$ contains no visible island on 13 points.  Indeed, for sake of contradiction, suppose $\hat{S}\subset \hat{H}_n$ is a visible island in $\hat{H}_n$ and $|\hat{S}| = 13$.
Since at most two members of a triplet can belong to $\hat{S}$, by the pigeonhole principle, there are 7 points $p_{i_1},\ldots, p_{i_7} \in H_n$ that are  parents of points in $\hat{S}$.  By setting $S = \{p_{i_1},\ldots, p_{i_7}\}$, Lemma \ref{lem} implies that there is a point $p_j \in H_n$ that lies in the interior of $\conv(S)$.  However, this implies that the triplet $q_j,u_j,v_j \in \conv(\hat{S})$, a contradiction.  $\hfill\square$

    \section{Concluding remarks}

Our initial goal was to find arbitrarily large visible islands in point sets with no four collinear members.  Unfortunately, Theorem \ref{main} shows that this is not possible.  However, the following conjecture remains open, which would imply the Big-line Big-clique Conjecture for $\ell = 4$ by applying an induction on $k$ and setting $n$ sufficiently large.  Given a finite point set $P$, the \emph{neighborhood} of $p\in P$ is the set of points in $P$ that are visible to $p$.

    \begin{conjecture}
Every $n$-element planar point set $P$ with no four collinear members, contains a point $p\in P$ such that its neighborhood contains an island of size $f(n)$, where $f(n)$ tends to infinity as $n$ tends to infinity.
    \end{conjecture}

    \medskip

        \noindent \textbf{Acknowledgement.} The first author would like to thank her sponsors for their generous support during the 2021 Undergraduate Summer Research Program at UCSD.

\end{document}